\documentclass[11pt,a4paper,fleqn]{article}

\usepackage{graphicx,natbib,amssymb,latexsym,lineno}

\newtheorem{prop}{Proposition}[section] 
\newenvironment{proof} {\noindent {\em \textbf{Proof}} } { \hfill $\Box$ \\ } 

\title {Properties of  Design-Based Functional Principal Components Analysis}

\author{Herv\'e Cardot, Mohamed Chaouch, Camelia Goga and Catherine Labru\`ere \\ 
Institut de Math\'ematiques de Bourgogne, UMR CNRS 5584 \\ Universit\'e de Bourgogne \\
 9 Avenue Alain Savary - B.P. 47870 \\
21078 DIJON Cedex - France}

\begin{document}

\maketitle

\begin{abstract}
This work aims at performing Functional Principal Components Analysis (FPCA) with Horvitz-Thompson estimators when the observations are curves  collected with survey sampling techniques.  One important motivation for this study is that FPCA is a dimension reduction tool which is the first step to develop  model assisted approaches that can take auxiliary information into account. FPCA relies on the estimation of the eigenelements of the covariance operator which can be seen as nonlinear functionals.  Adapting to our functional context the linearization technique based on the influence function developed by Deville (1999), we prove that these estimators are asymptotically design unbiased and consistent. Under mild assumptions, asymptotic variances are derived for the FPCA' estimators and  consistent estimators of them are proposed. Our approach is illustrated with a  simulation study and we check the good properties of the proposed estimators of the eigenelements as well as their  variance estimators obtained with the linearization approach.
\end{abstract}


\section{Introduction and notations}
Functional Data Analysis whose main purpose is to provide tools for describing and modeling  sets of curves is a topic of growing interest in the statistical community. The books by Ramsay and Silverman (2002, 2005) propose an interesting description of the available procedures dealing with functional observations whereas Ferraty and Vieu (2006) present  a completely nonparametric point of view. These functional approaches mainly rely on generalizing multivariate statistical procedures in functional spaces and have been proved useful in various domains such as  chemometrics (Hastie and Mallows, 1993), economy (Kneip and Utikal, 2001), climatology (Besse \textit{et al.} 2000), biology (Kirkpatrick and Heckman, 1989, Chiou \textit{et al.} 2003) or remote sensing (Cardot \textit{et al.}, 2003). These functional approaches are generally more appropriate than  longitudinal data models or time series analysis when there are for each curve many measurement points (Rice, 2004).

When dealing with functional data, the statistician generally  wants, in a first step, to represent as well as possible the sample of  curves in a  well chosen small dimension space in order to get a description of the functional data that allows interpretation. This objective can be achieved by performing a Functional Principal Components Analysis (FPCA) which provides  a small dimension space able to capture, in an optimal way according to a variance criterion, the main modes of variability of the data. These modes of variability are given  by considering, once the mean function has been subtracted off,  projections onto the space generated by  the eigenfunctions of the covariance operator associated to the largest eigenvalues. This technique is also known as the Karhunen-Loeve expansion in probability or Empirical Orthogonal Functions  (EOF) in climatology and numerous works have been published on this topic.  From a statistical perspective, the seminal paper by Deville (1974) introduces the functional framework whereas  Dauxois \textit{et al.} (1982) give asymptotic distributions. More recent works deal with smoothing or interpolation procedures (Castro \textit{et al.}, 1986, Besse and Ramsay, 1986, Cardot, 2000 or Benko \textit{et al.} 2009) as well as bootstrap properties (Kneip and Utikal, 2001) or sparse data (James \textit{et al.}, 2000).

The way data are collected is seldom taken into account in the literature and one generally supposes the data are independent realizations drawn from a  common functional probability distribution.
Even if this assumption can be supposed to be satisfied in most situations, there are some cases for which it will lead to estimation procedures that are not adapted to the sampling scheme. 
Design of experiments approaches have been studied by Cuevas \textit{et al.} (2003) but nothing has been done in the functional framework, as far as we know, from a survey sampling point of view whereas it can have some interest for practical applications. For instance, Dessertaine (2006) considers the estimation with time series procedures 
of electricity demand  at fine time scales with the observation of  individual electricity consumption curves. In this study, the data are  functions of time measured every ten minutes with more than 1000 time point observations and can be naturally thought as functional data.  Moreover, the individuals (\textit{e.g.} electricity meters) are selected according to balancing techniques (Deville and Till\'e, 2004) and consequently they do not  have the same probability to belong to  the sample. More generally, there are now data (data streams) produced automatically by large numbers of distributed sensors which generate huge amounts of data that can be seen as functional. The use of   sampling techniques to collect them proposed for instance in Chiky and H\'ebrail (2009) seems to be a relevant approach in such a framework allowing a trade off between limited storage capacities and accuracy of the data. In such  situations classical estimation procedures will lead to misleading interpretation of the FPCA since the mean and covariance structure of the data will not be estimated properly. 

We propose in this work estimators of  the FPCA when the curves are collected with survey sampling strategies. Let us note that Skinner \textit{et al.} (1986) have studied some properties of multivariate PCA in such a survey framework. Unfortunately, this work has received little attention in the statistical community. The functional framework is different since the eigenfunctions which exibit the main modes of variability of the data are also functions of time and can be naturally interpreted as modes of variability varying along time. FPCA can also be, by its dimension reduction properties, a useful tool  if one wants to use model-assisted approache (S\"arndal \textit{et al.}, 1992) that can take auxiliary information into account. Adapting for instance the single index model (Chiou \textit{et al.} 2003) or the additive model (M\"uller and Fang, 2008) on the principal components  scores  in this survey context would allow us to consider model assisted and small domain estimation in a functional context.

The paper is structured as follows. We first  define, in section 2, functional principal components analysis in a finite population setting. Then we propose estimators of the mean function and the covariance operator based on Horvitz-Thompson estimators. We also describe how this  dimension reduction tool can be of great interest for model assisted and small domain estimation when auxiliary information is available. Section 3 is devoted to the asymptotic properties. We show in section 3.1 that the FPCA' estimators are asymptotically design unbiased and consistent.  Section 3.2 provides approximations and consistent estimators of the variances of  FPCA' estimators with the help of perturbation theory (Kato, 1966) and the influence function (Deville, 1999). Campbell (1980) proposed, in a pioneer work, to use the influence function for estimating the variance of complex statistics and compared it with a jackknife variance estimator (see also Berger and Skinner, 2005). In such a functional context, we can not perform a first-order Taylor expansion  of the associated complex statistics but we can make a first-order von Mises (1947) expansion of the functional giving these complex statistics and obtain under broad assumptions  that the asymptotic variance of the complex statistics is equal to the variance of the Horvitz-Thompson estimator for the population total of some artificial variable $u_k$ constructed using the influence function technique (Deville, 1999). A jackknife variance estimator can be obtained by analogy with the Deville's linearization variance in which the analytic expression of $u_k$ is replaced by its numerical approximation (Davison and Hinkley, 1997). 
 Section 4 proposes a simulation study which shows the good behavior of our estimators for various sampling schemes as well as the ability of linearization techniques to give good approximations to their theoretical variances. The proofs are gathered in an Appendix.

\section{Survey framework and PCA}

\subsection{FPCA in a finite population setting}
Let us consider again the example of the estimation of the electricity demand presented in the introduction. If measures are taken every ten minutes during 24 hours, the consumption curve for one household $k$ belonging to the population is represented by the functional $Y_k(t)$ with $t$ being one of the 144 time measurements. In such a situation it is more convenient to consider that the observed trajectories are functions, instead of vectors of size 144, belonging to a function space that we suppose, from now on and without loss of generality, to be $L^2[0,1],$  the space of square integrable functions defined on the closed interval $[0,1].$ This space is equipped with the its inner product  $\langle \cdot, \cdot \rangle$ and norm $\| \cdot \|.$

Let us consider a finite population $U=\{1, \ldots, k, \ldots, N\}$ with size $N,$ not necessarily known,
and a functional variable $\mathcal{Y}$ defined for each element $k$ of the population $U$:
$Y_k=(Y_k(t))_{t \in [0,1]}$ belongs to the space $L^2[0,1].$  
Suppose first that we are looking for the function $\mu \in L^2[0,1]$ which is the closest to the population curves according to a quadratic loss criterion. The criterion
$ \sum_{k \in U} \left\| Y_k - \phi_0 \right\|^2 \ $ is clearly minimum  for $\phi_0=\frac{1}{N} \sum_{k \in U} Y_k,$ which is the mean population curve :
\begin{eqnarray}
\mu(t) & = &\frac{1}{N} \sum_{k \in U} Y_k(t), \quad t \in [0,1]\label{moyennefonct}
\end{eqnarray}

The curves $Y_k$ span a subspace of $L^2[0,1]$ whose dimension can be very large, at most $N.$
Going further, we would like now to obtain a subspace of $L^2[0,1]$ with dimension $q, $ as small as possible, that would allow to represent as well as possible the deviation of the population curves from their mean function $\mu.$   Considering an orthonormal basis $\phi_1, \phi_2, \ldots, \phi_q$ of this $q$ dimensional space, it is well known that the projection $P_q$ of the deviation of  the $Y_k$ from their mean function $\mu$ can be expressed as follows
$$
P_q(Y_k - \mu) = \sum_{j =1}^q \langle Y_k - \mu,\phi_j \rangle \phi_j .
$$
Considering again a quadratic  loss criterion, we would like to minimize the following quantity according to the set of orthonormal functions $\phi_1, \phi_2, \ldots, \phi_q,$
 \begin{eqnarray}
R_q (\phi_1, \phi_2, \ldots, \phi_q) & = & \frac{1}{N} \sum_{k=1}^{N} \left\| (Y_k -\mu) -  \left(\sum_{j =1}^q \langle Y_k - \mu,\phi_j \rangle \phi_j \right)  \right\|^2.
\label{criterionACP}
\end{eqnarray}
To get the solution of this optimization problem, we need to introduce more notations. Let us define the covariance operator, say $\Gamma,$ of the functions $Y_k, k \in U,$ as follows
\begin{eqnarray}
\Gamma &= &\frac{1}{N} \sum_{k \in U} \left(Y_k - \mu\right)\otimes \left(Y_k - \mu\right)\label{covariancefonct}
\end{eqnarray}
where the tensor product of two elements $a$ and $b$ of $L^2[0,1]$ is the rank one operator such that $a \otimes b(u) = \langle a, u \rangle b$ for all $u$ in $L^2[0,1].$
The operator $\Gamma$ is symmetric and non negative ($\langle \Gamma u, u \rangle \geq 0$). Its eigenvalues, which are positive and supposed to be sorted in decreasing order  $\lambda_1 \geq \lambda_2 \geq \cdots \geq \lambda_N \geq 0,$  satisfy  
\begin{equation}
\Gamma v_j(t) \ = \ \lambda_j \ v_j (t), \quad t \in [0,1], \quad j=1,\ldots, N, \label{valvecpropres}
\end{equation}
where the eigenfunctions $v_j, j=1, \ldots, q,$ form an orthonormal system in $L^2[0,1],$  \textit{i.e} $\langle v_j, v_{j'} \rangle = 1$ if $j=j'$ and zero otherwise. Going back to criterion (\ref{criterionACP}), one can express it as follows
\begin{eqnarray}
R_q (\phi_1, \phi_2, \ldots, \phi_q) &= & \frac{1}{N} \sum_{k \in U} \left\| Y_k - \mu \right\|^2 - \sum_{j=1}^q \langle \Gamma \phi_j, \phi_j \rangle
\end{eqnarray}
and we get by maximal properties of the eigenvalues (see Chatelin, 1983)  that the minimum of $R_q (\phi_1, \phi_2, \ldots, \phi_q)$ is  attained for $\phi_1=v_1, \ldots, \phi_q=v_q.$ Thus the optimal subspace of dimension $q,$ which is unique if  $\lambda_q > \lambda_{q+1},$  is the space generated by the $q$ eigenfunctions of $\Gamma$ associated to the $q$ largest eigenvalues.
Having these considerations in mind, we can build an expansion, which is similar to the Karhunen-Lo\`eve expansion or FPCA, that allows to get the best approximation in  a finite dimension space with dimension $q$ to the curves of the population making the key decomposition
 \begin{eqnarray}
Y_k(t) & = & \mu(t) + \sum_{j =1}^q \langle Y_k - \mu,v_j \rangle v_j(t) + R_{q,k}(t), \quad t \in [0,1]
\label{KLdecomp}
\end{eqnarray}
where $R_{q,k}(t)$ is the remainder term.
This means that  the space generated by the eigenfunctions $v_1, \cdots, v_q$ gives a representation of   the main modes of variation along time $t$ of the data around the mean $\mu.$ Moreover, the variance of the projection onto each $v_j$ is given by the eigenvalue
\begin{eqnarray*}
\lambda_j &= &\frac{1}{N} \sum_{k \in U} \langle Y_k - \mu,v_j \rangle^2 \ 
\end{eqnarray*}
since $\frac{1}{N} \sum_{k \in U} \langle Y_k - \mu,v_j \rangle=0.$

We aim, in the following, at estimating the mean function $\mu$ and the covariance operator $\Gamma$ in order to deduce estimators of the eigenelements $(\lambda_j, v_j)$ when the data are obtained with survey sampling procedures.  To this purpose, we express our parameters of interest as non-linear functions of finite population totals. Next, we substitute each total with its Horvitz-Thompson estimator described in the next section and finally, we obtain in section \ref{asymptoticprop} the asymptotic variance adapting  the influence function approach (Deville, 1999).

\noindent \textbf{Remark.} \textit{ In the space $L^2[0,1],$ we have in an equivalent way the following representation of the covariance operator
\begin{eqnarray}
\Gamma u(t) &=& \int_0^1 \gamma(s,t) u(s) \ ds\label{cov1}
\end{eqnarray}
where $\gamma(s,t)$ is the covariance function
\begin{eqnarray}
\gamma(s,t) &= &\frac{1}{N} \sum_{k \in U} \left(Y_k(t) - \mu(t)\right) \left(Y_k(s) - \mu(s)\right), \quad (s,t) \in [0,1]\times [0,1].\label{cov2}
\end{eqnarray}
Note also that if $H=\mathbb{R}^p$ then we get back to the classical definition of the principal components, the covariance operator being then the variance-covariance matrix, with size $p\times p,$ of the population vectors.}

\subsection{The Horvitz-Thompson Estimator}

Let us consider a sample $s$ of $n$ individuals, \textit{i.e.} a subset  $s\subset U,$  selected according to a probabilistic procedure $p(s)$ where $p$ is a probability distribution on the set of $2^N$ subsets of  $U.$ 
We denote by $\pi_k=\Pr(k\in s)$ for all  $k\in U$ the first order inclusion probabilities and by  $\pi_{kl}=\Pr(k \ \& \  l\in s)$ for all $k, l\in U$ with $\pi_{kk}=\pi_k, $ the second order inclusion probabilities.
 We suppose that all the individuals and all the pairs of individuals of the population have non null probabilities to be selected in the sample $s$, namely $\pi_k>0$ and  $\pi_{kl}>0$. We also suppose  that $\pi_k$ and $\pi_{kl}$ are not depending on $t\in [0,1]$. This means that once we have selected the  sample $s$ of individuals,  we observe $Y_k(t)$ for all $t\in [0,1]$ and all $k \in s.$
 Let us  start with the simplest case, the estimation of the finite population total of the $Y_k$ curves denoted by 
 $$
 t_Y=\sum_{k \in U} Y_k.
 $$ 
The Horvitz-Thompson (HT) estimator $\widehat t_{Y\pi}$ of  $t_Y$ is a function belonging to $L^2[0,1]$ defined as follows
\begin{eqnarray*}
\widehat t_{Y\pi}=\sum_{k \in s} \frac{Y_k}{\pi_k}=\sum_{k \in U} \frac{Y_k}{\pi_k}I_k\label{HTestim}
\end{eqnarray*}
where $I_k=\mathbf 1_{\{k\in s\}}$ is the sample membership indicator of element $k$ (S$\ddot{a}$rndal \textit{et al.}, 1992). Note that the variables $I_k$ are random with $Pr(I_k=1)=\pi_k$ whereas the curves $Y_k$ are considered as  fixed with respect to the sampling design $p(s)$. So, the HT estimator $\widehat t_{Y\pi}$ is $p$-unbiased, namely
$$
E_p(\widehat t_{Y\pi})=t_Y
$$
where $E_p(\cdot)$ is the expectation with respect to the sampling design.

The variance operator of $\widehat t_{Y\pi}$ calculated with respect to $p(s)$ is the HT variance 
\begin{eqnarray}
\mbox{V}_p(\widehat t_{Y\pi})=\sum_U\sum_U \Delta_{kl}\frac{Y_k}{\pi_k}\otimes\frac{Y_l}{\pi_l}\label{HTvar}
\end{eqnarray}
and it is estimated $p$-unbiasedly by $\displaystyle\widehat{\mbox{V}}_p(\widehat t_{Y\pi})=\sum_s\sum_s \frac{\Delta_{kl}}{\pi_{kl}}\frac{Y_k}{\pi_k}\otimes\frac{Y_l}{\pi_l}$
with the notation $\Delta_{kl}=\pi_{kl}-\pi_k\pi_l$ if $k \neq l$ and $\Delta_{kk}=\pi_k(1-\pi_k).$ One may obtain equivalent integral representations of $\mbox{V}_p(\widehat t_{Y\pi})$ and $\widehat{\mbox{V}}_p(\widehat t_{Y\pi})$ similar as in equations (\ref{cov1}) and (\ref{cov2}). 

\noindent\textbf{Example}: \textit{Let us select a sample of $n$ curves $Y_k$ according to a simple random sample without replacement (SI)  from $U$. We have $\widehat t_{Y\pi}=(N/n)\sum_{k \in s} Y_k$ with variance $\mbox{V}_{SI}(\widehat t_{Y\pi})=N^2\frac{1-f}{n}S^2_{YU}$ for $f=n/N$ and $S^2_{YU}=\frac{1}{N-1}\sum_U(Y_k-\mu)\otimes(Y_k-\mu)$ the population variance. The variance estimator is given by $\widehat{\mbox{V}}_{SI}(\widehat t_{Y\pi})=N^2\frac{1-f}{n}S^2_{Ys}$ with $S^2_{Ys}=\frac{1}{n-1}\sum_s(Y_k-\mu_s)\otimes(Y_k-\mu_s)$ and $\mu_s=\frac{1}{n}\sum_s Y_k$}.

\subsection{Substitution Estimator for Nonlinear Parameters}

 Consider now  the estimation of one of the following parameters: $\mu, $ $\Gamma$ and the eigenelements $\lambda_j$ and $v_j$ given by  (\ref{moyennefonct}), (\ref{covariancefonct}) and (\ref{valvecpropres}). When the population size is unknown, we deal with  nonlinear functions of population totals.
To estimate these parameters, we substitute each total by its Horvitz-Thompson estimator as described in the above section. We obtain  complex statistics whose variances are no longer calculated using formula (\ref{HTvar}). Besides the nonlinearity feature, we have to cope now with the fact that $\mathcal Y$ is a functional variable which makes the variance estimation issue more difficult. In order to overcome this, we adapt  the linearization technique based on the influence function  introduced by Deville (1999)  to  the functional framework. This approach is based on the fact that each finite population total may be written as a functional  depending on a finite and discrete measure $M$ and as a consequence, the population parameter of interest can be written as a functional $T(M). $ We derive the Horvitz-Thompson estimator $\widehat M$ of $M$ and estimators or our parameters are obtained by pluging-in $\widehat M$ in the expression of $T,$ namely $T(\widehat M).$\\
Let us introduce now  the discrete measure $M$ defined on $L^2[0,1]$ as follows 
\begin{eqnarray*}
M & = & \sum_{k \in U} \delta_{Y_k}
\end{eqnarray*}
where $\delta_{Y_k}$ is the Dirac function taking value 1 if $Y=Y_k$ and zero otherwise.
The following  parameters of interest can be defined as functionals of $M$:
\begin{eqnarray*}
N & = & \int dM \quad \mbox{and}\quad\mu  =   \frac{\displaystyle \int \mathcal{Y} dM}{\displaystyle \int dM} \\
\Gamma &=&  \frac{\displaystyle \int \left( \mathcal{Y}-\mu \right) \otimes   \left( \mathcal{Y}-\mu \right) dM}{\displaystyle \int dM}
\end{eqnarray*}
\noindent and the eigenelements, given by (\ref{valvecpropres}), are implicit functionals $T$ of $M$.\\
\noindent The measure $M$ is estimated by the Horvitz-Thompson estimator $\widehat{M}$ associating the weight $1/\pi_k$ for each $Y_k$ with $k\in s$ and zero otherwise,
\begin{eqnarray*}
\widehat{M} & = & \sum_{k \in U} \frac{\delta_{Y_k}}{\pi_k} \ I_k
\end{eqnarray*}
and $T(M)$ is then estimated by $T(\widehat M)$ also called the \textit{substitution estimator}. For example, the substitution estimators for $\mu$ and  $\Gamma$ are 
\begin{eqnarray}
\widehat{\mu}  & = & \frac{1}{\widehat{N}} \sum_{k \in s} \frac{Y_k}{\pi_k} \label{muest} \\
\widehat{\Gamma} & = & \frac{1}{\widehat{N}} \sum_{k\in s} \frac{Y_k\otimes Y_k}{\pi_k} -  \widehat{\mu}\otimes\widehat{\mu}
\label{gammaest}
\end{eqnarray}
where the size $N$ of the population is estimated by  $\widehat{N}  =  \displaystyle \sum_{k \in s} \frac{1}{\pi_k}.$
Then estimators of the eigenfunctions $\{\widehat{v}_j, j=1 , \ldots q\}$ associated to the $q$ largest  eigenvalues  $\{\widehat{\lambda}_j, j=1 , \ldots q\}$ are obtained readily by the  eigen-analysis of the estimated covariance operator $\widehat{\Gamma}.$

\textbf{Remark.} \textit{In practice we do not observe the whole curves  but generally discretized versions at $m$ design points $0 \leq t_1 < t_2 < \cdots < t_m \leq 1$ that we suppose to be the same for all the curves. Quadradure rules are often employed in order to get numerical approximations to integrals and inner product by summations : for each $u$ in $L^2[0,1]$ we get  an accurate discrete approximation to the integral
$$
\int_0^1 u(t) dt \approx \sum_{\ell=1}^m w_\ell \ u(t_\ell)
$$
provided the number of design points $p$ is large enough and the grid is sufficiently  fine. When the discretization points vary from one curve to another basis functions approaches are generally employed in order to smooth and to decompose the signals in a common functional space (see e.g. Ramsay and Silverman, 2005).}

\subsection{Some comments on the interest of FPCA in survey sampling}

As seen in equation (\ref{KLdecomp}), the FPCA allows to get a finite and generally small dimension space that is able to reconstruct rather well the curves of the population. Indeed the principal components scores $\langle Y_k - \mu, v_j \rangle,$ for $j=1, \ldots, q,$ are indicators of the deviation of curve $Y_k$ from its mean function $\mu.$  When auxiliary variables that influence significantly the shape of the population curves are known for each element of the population it would certainly be of great interest to consider models that could explain the individual fluctuations of the principal components scores. This could be useful for instance for improving the total curve estimation and small domain estimation. Indeed, suppose we have a set of $p$ real covariates, $x_{1}, \ldots, x_{p}$ and a function $f_j$ (to be estimated) such that
$$
\xi : \quad \langle Y_k - \mu, v_j \rangle = f_j(x_{k1}, \ldots, x_{kp}) + \epsilon_{jk}
$$
where $\epsilon_{jk}$ is supposed to be a random noise for $k \in U$ and $j=1, \ldots, q.$ 
Then, having built estimators $\widehat{f}_j$ of the functions $f_j$ using the model $\xi$ and the sampling design $p(\cdot),$ the total curve could be estimated considering a model assisted approach
$$
\widehat{t}_Y (t) = \sum_{k \in s} \frac{Y_k (t)}{\pi_k} - \left(\sum_{k \in s} \frac{\widehat Y_k (t)}{\pi_k}-\sum_U\widehat Y_k(t)  \right), \quad t \in [0,1].
$$
where the predicted $Y$'s values are given by 
$$\widehat Y_k(t)=\widehat{\mu}(t) + \sum_{j=1}^q \widehat f_j(x_{k1}, \ldots, x_{kp})\widehat v_j (t).$$
Chiou  \textit{et al.} (2003) proposed nonparametric estimators of function $f_j$ based on single index models that could explain the principal components scores thanks to real covariates whereas M\"uller and Yao (2008) consider additive models.

Going back to the motivating example of individual electricity consumption curves, it is clear that the temperature, the past consumption, or the  surface of the household can be made available for the all population and are certainly correlated with the shape of individual curves. Thus building statistical models that explain the variations of the principal components scores should be helpful to provide better estimators of the total consumption curve as well as estimation of total curves for small domains.
This issue which is according to us of great interest deserves further investigations that are beyond the scope of this paper.

\section{Asymptotic Properties}\label{asymptoticprop}
We give in this section asymptotic properties of our estimators $\hat \mu,$ $\hat \Gamma$ and $\hat\lambda_j, \hat v_j.$ Nevertheless, the approach we propose in the following is general and can be useful for estimating other non-linear functions of totals.\\
Let us consider the superpopulation asymptotic framework introduced by Isaki and Fuller (1982) which supposes that the population and the sample sizes tend to infinity.  Let $U_{\mathbb{N}}$ be a population  with infinite (denumerable) number of individuals and consider a sequence of nested sub populations such that $U_1\subset \cdots  \subset U_{\nu-1} \subset U_{\nu} \subset U_{\nu+1} \subset \cdots \subset U_{\mathbb{N}}$ of sizes $N_1<N_2<\ldots<N_{\nu}<\ldots$.
Consider then a  sequence of  samples $s_{\nu}$ of size $n_{\nu}$  drawn from $U_{\nu}$ according to the fixed-size sampling designs $p_{\nu}(s_{\nu})$ and denote by $\pi_{k\nu}$ and $\pi_{kl\nu}$ their  first and second order inclusion probabilities . Note that the sequence of sub populations is an increasing nested one while the sample sequence is not. For sake of simplicity, we will drop the subscript $\nu$ in the following.

We assume that the following assumptions are satisfied :
\begin{itemize}
\item[\textbf{(A1)}] $\quad \displaystyle \sup_{k \in U} \left\| Y_k \right\| \leq C <\infty$,
\item[\textbf{(A2)}] $\displaystyle \lim_{N\rightarrow\infty}\frac{n}{N}=\pi\in(0,1) $,
\item[\textbf{(A3)}] $\displaystyle \min_{k \in U_N} \pi_k\geq \lambda>0 \ , \quad \min_{k\neq l} \pi_{kl}\geq \lambda*>0$ 
and  $\displaystyle \overline{\lim}_{N\rightarrow\infty}n \max_{k\neq l} |\pi_{kl}-\pi_k\pi_l| <\infty.$
\end{itemize}
Hypothesis (A1) is rather classical in functional data analysis. Note that  it does not imply that the curves $Y_k(t)$ are uniformly bounded in  $k$ and $t \in [0,1].$ 
 Hypotheses (A2) and (A3)  are checked for usual sampling plans (Robinson and S\"arndal, 1983, Breidt and Opsomer, 2000).

 \subsection{ADU-ness and Consistency of Estimators}
The substitution estimators of  $\mu$ and $\Gamma$ defined in (\ref{muest}) and (\ref{gammaest}), as well as $\hat\lambda_j$ and $\hat v_j$, are no longer $p$-unbiased. Nevertheless, we show in the next that, in large samples, they are \textit{ asymptotically design unbiased} (ADU) and \textit{consistent}.\\
An estimator $\widehat\Phi$  of $\Phi$  is said to be \textit{ asymptotically design unbiased} (ADU) if 
$$
\lim_{N\rightarrow\infty}\left( E_p(\widehat\Phi)-\Phi\right)=0.
$$
We say that  $\widehat\Phi$  satisfies $\left(\widehat\Phi - \Phi\right) = O_p(u_n)$ for a sequence $u_n$ of positive numbers if there is a constant $C$ such that  for any $\varepsilon>0,$ $\Pr \left( \left| \widehat\Phi - \Phi \right| \geq C  u_n \right) \leq \epsilon.$ The estimator is \textit{consistent} if one can find a sequence $u_n$ tending to zero as $n$ tends to infinity such as $\widehat\Phi - \Phi  = O_p(u_n).$ \\
 Let us also introduce the Hilbert-Schmidt norm, denoted by $\left\| \cdot \right\|_2$  for operators mapping $L^2[0,1]$  to $L^2[0,1].$ It is induced by the inner product between two operators $\Gamma$ and $\Delta$ defined by
$\langle  \Gamma, \Delta\rangle_2 =  \sum_{\ell=1}^{\infty} \langle \Gamma e_\ell, \Delta e_\ell \rangle$ for any orthonormal basis $(e_\ell)_{\ell \geq 1}$ of $L^2[0,1].$ In particular, we have that 
$ \left\| \Gamma \right\|_{2}^2 = \sum_{\ell=1}^{\infty} \langle \Gamma e_\ell, \Gamma e_\ell \rangle = \sum_{j\geq 1} \lambda_j^2.$ 
 
\begin{prop} Under hypotheses (A1), (A2) and (A3),
\begin{eqnarray*}
E_p \left(\frac{N - \widehat{N}}{N}\right)^2 & =& O( n^{-1}), \\
E_p \left\| \mu - \widehat{\mu} \right\|^2 & =& O( n^{-1}), \\
E_p \left\| \Gamma - \widehat{\Gamma} \right\|_{2}^2 & =& O( n^{-1}).
\end{eqnarray*}
If we suppose that the non null eigenvalues are distinct, we also have,
\begin{eqnarray*}
E_p \left(\sup_j \left| \lambda_j - \widehat{\lambda_j} \right|\right)^2  & =& O( n^{-1}), 
\end{eqnarray*}
and for each fixed j,
\begin{eqnarray*}
E_p \left\| v_j - \widehat{v_j} \right\|^2 & =& O( n^{-1}) . 
\end{eqnarray*}
As a consequence, the above estimators are ADU and consistent.
\label{distmesure}
\end{prop}
\vspace{-0.3cm}
The proof is given in the Appendix. 

\subsection{Variance Approximation and Estimation}
Let us now define, when it exists, the influence function of a functional $T$ at point $\mathcal Y\in L^2[0,1]$ say $IT(M,\mathcal Y),$ as follows
 \begin{eqnarray*}
 IT(M, \mathcal Y) &= &\displaystyle \lim_{h \rightarrow 0} \frac{T(M + h \delta_\mathcal Y) - T(M)}{h} 
 \end{eqnarray*}
 where $\delta_\mathcal Y$ is the Dirac function at $\mathcal Y.$ Note that this  is not exactly the usual definition of the influence function (see \textit{e.g.} Hampel, 1974 or Serfling, 1980) and it has been adapted to the survey sampling framework by  Deville (1999). We define the \textit{linearized variables} $u_k$, $k\in U$ as the influence function of $T$ at $M$ and $\mathcal Y=Y_k$, namely
\begin{eqnarray*}
 u_k=IT(M,Y_k).\label{linvar}
 \end{eqnarray*}
 Note that the linearized variables depend on $Y_k$ for all $k\in U$ and as a consequence, they are all unknown.\\

\noindent We can give  a first order von Mises expansion of our functional $T,$ 
\begin{eqnarray}
T(\widehat{M})  & = &  T(M) + \sum_{k \in U} IT\left(M,Y_k\right) \left(\frac{I_k}{\pi_k} - 1\right) + R_T\\\nonumber
 & & = T(M) + \sum_{k \in U} u_k \left(\frac{I_k}{\pi_k} - 1\right) + R_T\label{IF:L2}
\end{eqnarray}
for the reminder $R_T$ and the linearized variable $u_k=IT(M,Y_k).$ The above expansion tells us, under regularity conditions, that the asymptotic variance of the estimator $T(\widehat{M})$ is the variance of  the HT estimator of the population total of  $IT\left(M,Y_k\right)$ provided the remainder term $R_T$ is negligible.

\noindent Before handling the remainder term, let us first calculate the influence function for our parameters of interest.
\begin{prop} Under assumption (A1), we get that  the influence functions of $\mu$ and $\Gamma$ exist and 
  \begin{eqnarray}
 I \mu (M, Y_k)  & =  &\frac{1}{N}(Y_k - \mu)\label{lineariseemu} \\
 I \Gamma (M,Y_k) &= &\frac{1}{N} \left( (Y_k-\mu)\otimes (Y_k-\mu) - \Gamma \right).\nonumber 
\end{eqnarray}
If moreover, the non null eigenvalues of  $\Gamma$ are distinct,  then
\begin{eqnarray}
I \lambda_j (M,Y_k) &= & \frac{1}{N} \left( \langle Y_k - \mu, v_j \rangle^2 - \lambda_j \right) \label{lineariseelambda}\\
I v_j (M,Y_k) &= & \frac{1}{N} \left( \sum_{\ell \neq j} \frac{\langle Y_k - \mu, v_j \rangle \langle Y_k - \mu, v_\ell \rangle}{\lambda_j - \lambda_\ell} v_\ell\right).\label{lineariseev}
\end{eqnarray} 
\label{IFacpf}
\end{prop}

 The proof is given in the Appendix. Let us remark that the influence functions of the eigenelements are similar to those found in the multivariate framework for classical PCA (Croux and Ruiz-Gazen, 2005). \\
We are now able to state that the remainder term $R_T$ defined in equation (\ref{IF:L2}) is negligible and  that the linearization approach can be used to get the  asymptotic variance of our substitution estimators. \\
 Let us suppose the supplementary assumption:
\begin{itemize}
\item[\textbf{(A4)}]: The Horvitz-Thompson estimator $\sum_s\frac{u_k}{\pi_k}$ satisfies a Central Limit Theorem  for the linearized variables  $u_k$  given by proposition \ref{IFacpf}.
\end{itemize}
 This assumption is satisfied for classical sampling designs and real quantities $u_k$ (see \textit{e.g} Chen and Rao, 2007, and references therein). The case of functional quantities deserves further investigations. Preliminary results can be found in Cardot and Josserand (2009).
  
 \begin{prop} Suppose the hypotheses (A1), (A2) and (A3) are true. Consider the functional $T$ giving the parameters of interest defined in (\ref{moyennefonct}), (\ref{covariancefonct}) and (\ref{valvecpropres}). We suppose that the non null eigenvalues are distinct. Then $
R_T  = o_p( n^{-1/2})$ and
\begin{eqnarray*}
T(\widehat{M}) - T(M)  =  \sum_{k \in U}   u_k \left(\frac{I_k}{\pi_k} - 1\right) + o_p(n^{-1/2})\label{var:IF}
\end{eqnarray*}
where the  $u_k$ are linearized variable of $T$ calculated in proposition \ref{IFacpf}.

If (A4) is also true, the asymptotic variance of $\widehat \mu$, resp. of $\widehat v_j$, is equal to the variance operator of the HT estimator  
$\displaystyle\sum_s\frac{u_k}{\pi_k}$ with $u_k$ given by (\ref{lineariseemu}), resp. by (\ref{lineariseev}), and its expression is given by 
\begin{eqnarray}
AV_p(T(\widehat M))=\sum_U\sum_U\Delta_{kl}\frac{u_k}{\pi_k}\otimes\frac{u_l}{\pi_l}
\label{varestimmuvect}
\end{eqnarray} 
The asymptotic variance of $\widehat\lambda_j$ is  
\begin{eqnarray}
AV_p(\widehat\lambda_j)=\displaystyle \sum_U\sum_U\Delta_{kl}\frac{u_k}{\pi_k}\frac{u_l}{\pi_l}\label{varestimlambda}
\end{eqnarray}  with $u_k$ given by (\ref{lineariseelambda}).
\label{Varapprox}
\end{prop}
The proof is given in the Appendix.
As one can notice, the asymptotic variances given in Proposition \ref{Varapprox} are unknown since the double sums are considered on the whole population $U$ and we have only a subset of it and secondly, the linearized variables $u_k$ are not known. As a consequence, we propose to estimate (\ref{varestimmuvect}) and (\ref{varestimlambda}) by the HT variance estimators replacing the linearized variables by their estimations. 
 In the case of $\widehat \mu$, $\widehat\lambda_j$ and  $\widehat v_j$, we obtain the following variance estimators: 

 \begin{eqnarray*}
\widehat{V}_p (\widehat{\mu}) & = & \frac{1}{\widehat{N}^2} \sum_{k \in s} \sum_{\ell \in s} \frac{1}{\pi_{k \ell}} \frac{\Delta_{k \ell}}{\pi_k \pi_\ell} \left( Y_k -\widehat{\mu} \right) \otimes \left( Y_\ell -\widehat{\mu} \right) \\
\widehat{V}_p \left(\widehat{\lambda}_j\right) & = & \frac{1}{\widehat{N}^2} \sum_{k \in s} \sum_{\ell \in s} \frac{1}{\pi_{k \ell}} \frac{\Delta_{k \ell}}{\pi_k \pi_\ell}
\left(  \langle Y_k - \widehat{\mu}, \widehat{v}_j \rangle^2 - \widehat{\lambda}_j \right)\left(  \langle Y_\ell - \widehat{\mu}, \widehat{v}_j \rangle^2 - \widehat{\lambda}_j \right)  \\
\widehat{V}_p \left( \widehat{v}_j \right) & = &  \sum_{k \in s} \sum_{\ell \in s}  \frac{1}{\pi_{k \ell}} \frac{\Delta_{k \ell}}{\pi_k \pi_\ell} \widehat{Iv}_j(M,Y_k) \otimes \widehat{Iv}_j(M,Y_\ell),
\end{eqnarray*}
with $\widehat{Iv}_j(M,Y_\ell) = \displaystyle \frac{1}{\widehat{N}} \left( \sum_{\ell \neq j} \frac{\langle Y_k - \widehat\mu, \widehat{v}_j \rangle \langle Y_k - \widehat\mu, \widehat{v}_\ell \rangle}{\widehat\lambda_j - \widehat\lambda_\ell} \widehat{v}_\ell\right).$

In order to prove  that these variance estimators are consistent we need to introduce additional assumptions involving higher order inclusion probabilities. 
\begin{itemize}
\item[\textbf{(A5)}] : Denote by $D_{t,N}$ the set of all distinct $t$ tuples $(i_1,i_2, \ldots, i_t)$ from $U.$ We suppose that
\begin{eqnarray*}
\lim_{N \rightarrow \infty} n^2 \max_{(i_1,i_2,i_3,i_4) \in D_{4,N}}  \left| E_p \left[ (I_{i_1} - \pi_{i_1}) (I_{i_2} - \pi_{i_2})(I_{i_3} - \pi_{i_3})(I_{i_4} - \pi_{i_4})\right] \right|  &< & \infty \\
\lim_{N \rightarrow \infty}  \max_{(i_1,i_2,i_3,i_4) \in D_{4,N}}  \left| E_p \left[ (I_{i_1}I_{i_2} - \pi_{i_1i_2}) (I_{i_3}I_{i_4} - \pi_{i_3i_4})\right] \right| & =& 0  \\
\lim \sup_{N \rightarrow \infty} n \max_{(i_1,i_2,i_3) \in D_{3,N}} \left| E_p \left[ (I_{i_1} - \pi_{i_1})^2 (I_{i_2} - \pi_{i_2})(I_{i_3} - \pi_{i_3})\right] \right|  &< & \infty 
\end{eqnarray*}
\end{itemize}
Hypothesis (A5) is a technical assumption that is similar to assumption A7 in Breidt and Opsomer (2000). These authors  explain in an interesting discussion that this set of assumptions holds for instance for simple random sampling without replacement (SRSWR) and stratified sampling.

\begin{prop}\label{estimvarconv} Under assumptions (A1)--(A5), we have that
\begin{eqnarray*}
E_p\left\| AV_p(\widehat{\mu})-\widehat{V}_p(\widehat{\mu}) \right\|_{2} &=& o \left(\frac{1}{n}\right) \\
E_p\left| AV_p(\widehat{\lambda}_j)-\widehat{V}_p(\widehat{\lambda}_j) \right| &=& o \left(\frac{1}{n}\right) 
\end{eqnarray*}
If moreover $\Gamma$ is a finite rank operator whose rank does not depend on $N$  then
\begin{eqnarray*}
\left\| AV_p(\widehat{v}_j)-\widehat{V}_p(\widehat{v}_j) \right\|_{2} &=& o_p \left(\frac{1}{n}\right) 
\end{eqnarray*} 
for $j=1, \ldots, q.$ 

\end{prop}
The proof is given in the Appendix. This theorem implies that  variance estimators for the mean function, the eigenvalues and the first $q$ eigenfunctions  are asymptotically design unbiased and consistent. 
Note that the hypothesis that $\Gamma$ is a finite rank operator 
is a technical assumption that is needed in the proof for the eigenfunctions in order to counterbalance the fact that eigenfunction estimators are getting poorer as $j$ increases. Note that with finite populations, operator $\Gamma$ is always a finite rank operator and its rank is at most $N,$ the population size. We probably  could  assume, at the expense of  more complicated proofs, that the rank of $\Gamma$ tends to infinity as $N$ increases. Allowing then $q$ to tend to infinity with the sample size with a rate depending on the shape of the eigenvalues should lead to the same variance approximation results for the eigenvectors.

\section{A simulation study}
We check now with a simulation study that we get accurate estimations to the eigenelements  even for moderate sample sizes as well as good approximation to their variance for simple random sampling without replacement (SRSWR) and stratified sampling.  In our simulations all functional variables are discretized in $m=100$ equispaced points in the interval $[0,1].$  Riemann approximations to the integrals are employed to deal with the discretization effects. 

We consider a random variable $Y$ following a Brownian motion with mean function $\mu(t)=\cos(4\pi t), t\in [0,1]$ and covariance function $cov(s,t)=\min(s,t)$. We make $N=10000$ replications of $Y$. We construct then two  strata $U_1$ and $U_2$ of different variances by multiplying the $N_1=7000$ first replications of $Y$ by $\sigma_1=2$ and the $N_2=3000$ other replications by $\sigma_2=4$. Our population $U$ is the union of these two strata.

\begin{figure}[h!]
\begin{center}
\includegraphics[height=13cm,width=16cm]{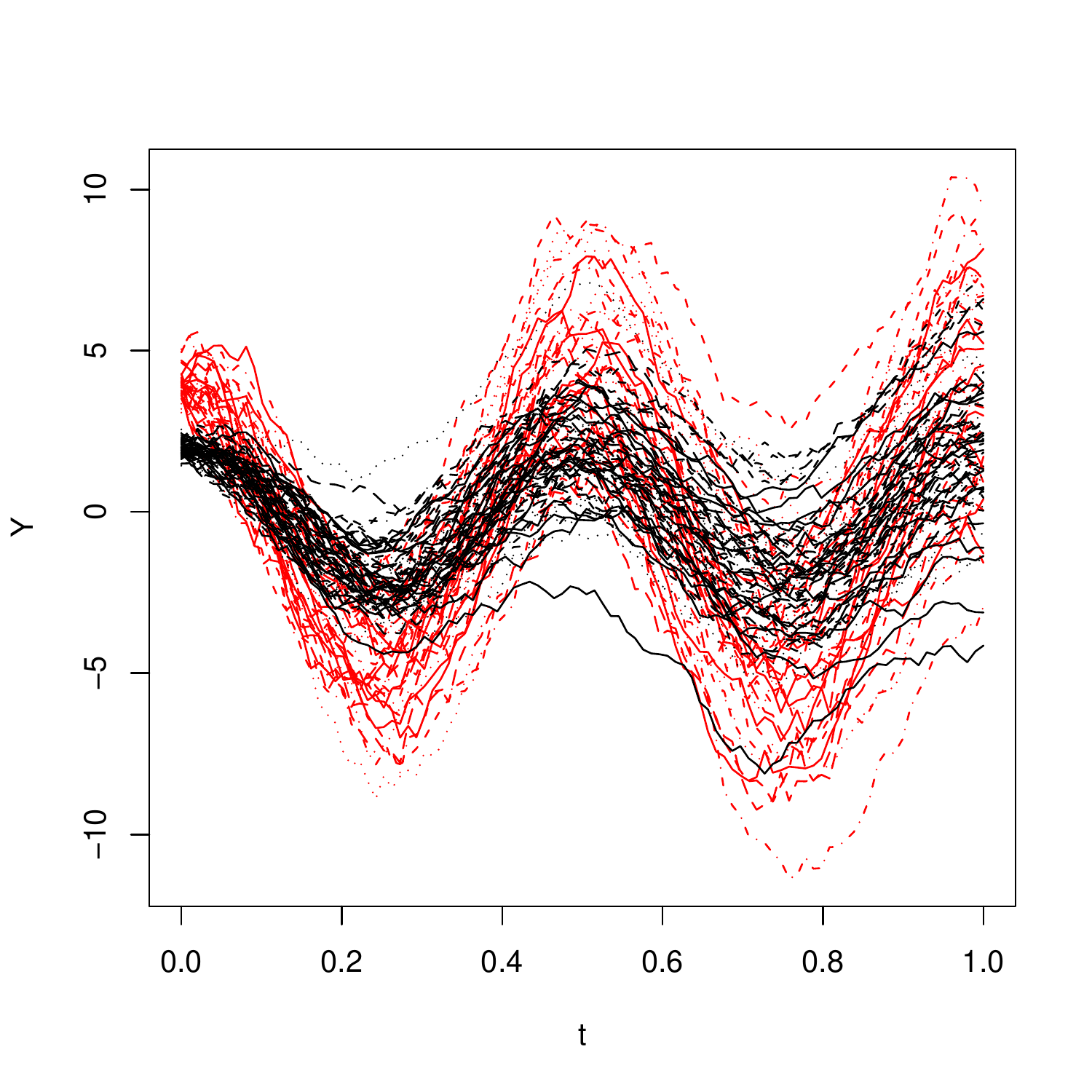} 
\end{center}

\caption{A stratified sample of $n=100$ curves}
\label{figpop}
\end{figure}

To evaluate our estimation procedures we make 500 replications of the following experiment.
We draw samples according to two different sampling designs (SRSWR  and stratified) and  consider two different sample sizes $n=100$ and $n=1000$. Each stratified sample is built by  drawing independently two  SRSWR of  sizes $n_1$ in stata $U_1$ and $n_2=n-n_1$ in  strata $U_2.$ The sample  sizes are chosen  to take into account the different variances in the strata: $$\frac{n_1}{n}=\frac{N_1}{N}\frac{\sigma_1}{\frac{ N_1 \sigma_1+ N_2 \sigma_2}{N}}\  ,\  \frac{n_2}{n}=\frac{N_2}{N}\frac{\sigma_2}{\frac{ N_1 \sigma_1+ N_2 \sigma_2}{N}}$$ in analogy with univariate stratified sampling with optimal allocation (S\"arndal \textit{et al.}, 1992). A stratified sample $s$ of size $n=100$ trajectories is drawn in Figure \ref{figpop}.

\begin{figure}[h!]
\begin{center}
\includegraphics[height=13cm,width=16cm]{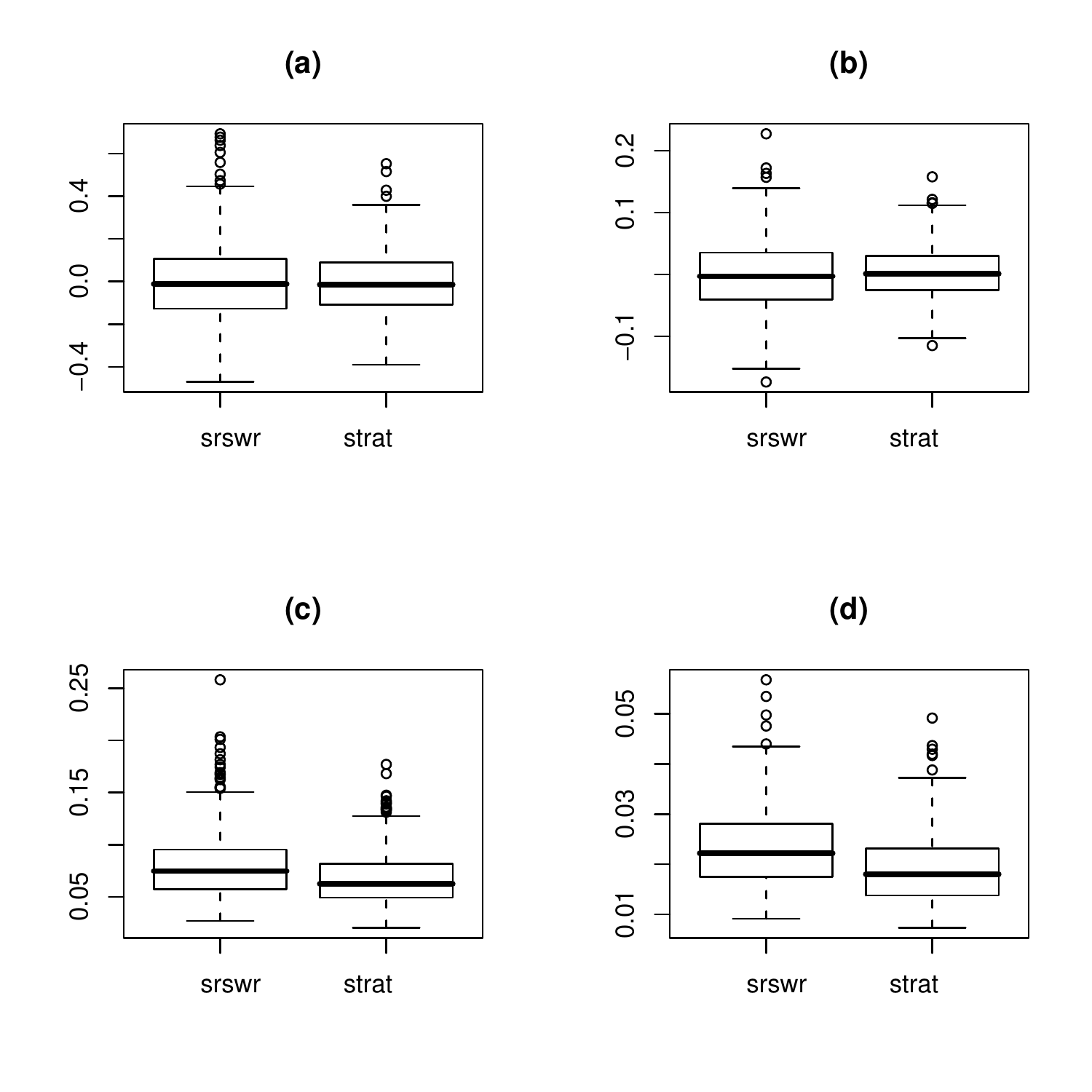} 
\end{center}
\caption{Estimation errors for two different sampling strategies (SRSWR and stratified sampling). First eigenvalue with $n=100.$ (a) and $n=1000$ (b). First eigenvector with $n=100.$ (c) and $n=1000$ (d).}
\label{figest}
\end{figure}

Estimation errors  for the first eigenvalue and the first eigenvector are evaluated  by considering the following loss criterions 
$\frac{\lambda_1-\widehat{\lambda_1}}{\lambda_1}$ and 
$\frac{\|v_1-\widehat{v_1}\|}{\|v_1\|}$ (Euclidean norm) among our 500 replications of the experiments.
The approximations turn out to be effective as seen in Figure \ref{figest}. For example for both sampling strategies the first eigenvector approximation has a median error lower than $3\%$ for a sample size $n=1000$. It also appears that the stratified sampling gives better estimations than the SRSWR sampling.

\begin{figure}[h!]
\begin{center}
\includegraphics[height=13cm,width=16cm]{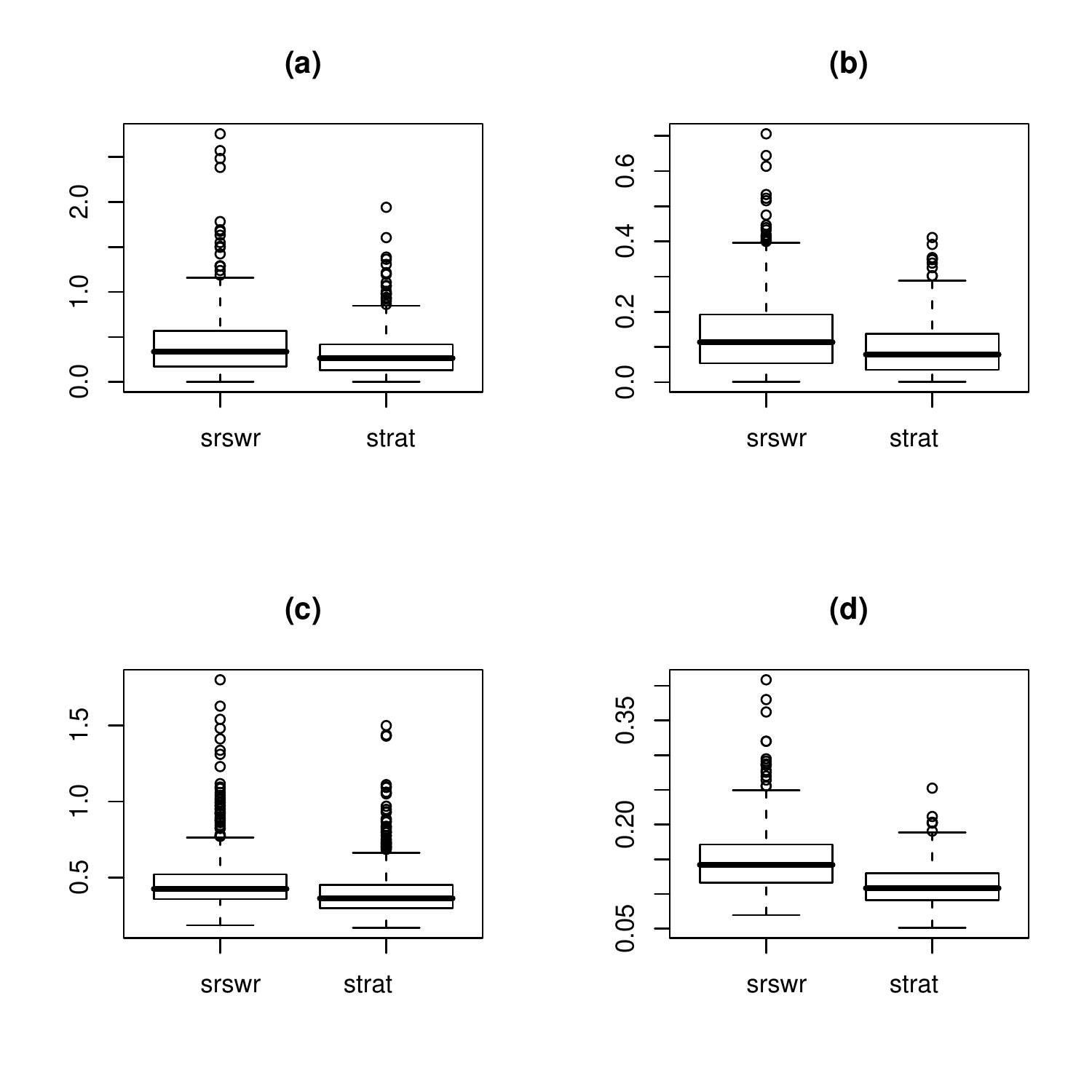} 
\end{center}

\caption{Estimation errors in the variance approximation for two different sampling strategies (SRSWR and stratified sampling). First eigenvalue with $n=100.$ (a) and $n=1000$ (b). First eigenvector with $n=100.$ (c) and $n=1000$ (d).}
\label{figvar}
\end{figure}

Let us look now at the variance of our estimators. Tables \ref{varvalp} and \ref{varvectp} give three variance (resp. euclidean norm of variance) approximations to the estimator of  respectively the first eigenvalue and the first eigenvector.  The first variance approximation to these estimators is their  empirical variance and are denoted by $Var(\widehat{\lambda_1})$ and by $Var(\widehat{v_1})$ , the second one is the  asymptotic variance  denoted by $AV(\widehat{\lambda_1})$ and by $AV(\widehat{v_1})$  whereas the third one is a $[25\%,75\%]$ confidence interval obtained by estimating  the asymptotic variance using the HT variance estimator respectively denoted by $\widehat{V}_p \left(\widehat{\lambda}_1\right)$ and $\widehat{V}_p \left(\widehat{v}_1\right)$ .  Errors (see Figure \ref{figvar})  in approximating the variance of the estimators by the  linearization approach are evaluated  by considering the following criterions: $\left |\frac{Var(\widehat{\lambda_1})-\widehat{V}_p \left(\widehat{\lambda}_1\right)}{Var(\widehat{\lambda_1})}\right |$ and 
 $\frac{\| Var(\widehat{v_1})-\widehat{V}_p \left(\widehat{v}_1\right)\|}{\|Var(\widehat{v_1}) \|  }.$
 
As a conclusion, we first note with this simulation study that HT estimators of the covariance structure of functional observations are accurate enough to derive good estimators of the FPCA. Secondly, linear approximations by the influence function give reasonable estimation of the variance of the eigenelements for small sample sizes and accurate estimations as far as $n$ gets larger (n=1000). We also notice that the variance of the estimators obtained by stratified sampling turns out to be  smaller than with SRSWR sampling.

%%%%%%%%%%%%%%%%%%%%%%%%%%%%%%%%%%%%%%%%%

%%%%%%%%%%%%%%%%%%%%%%%%%%%%%%

%%%%%%%%%%%%%%%%%%%%%%%%%%%%%%%%%%%%%%%%%%%%%%%%%%%%%%%%%%%%%%%%%%%%%%%%%%%%%%%%%%%%%%%%%
\begin{table}[h!]
	\centering
		\begin{tabular}{|c|c|c||c|c|}
		\hline
			& \multicolumn{2}{|c||}{n=100} & \multicolumn{2}{c|}{n=1000}\\
			\hline
			&SRSWR&stratified&SRSWR&stratified\\
			\hline
			$Var(\widehat{\lambda_1})$&0.314&0.223&0.0317&0.0189\\
			\hline
			 $AV(\widehat{\lambda_1})$&0.340&0.209&0.0309&0.0183\\
			\hline
			 $\widehat{V}_p \left(\widehat{\lambda}_1\right)$ & [0.208;0.430]&[0.155;0.257]&[0.027;0.034]&	[0.0169;0.0195]\\
			\hline
		\end{tabular}
		\caption{Variance approximation of the first eigenvalue estimator. }
		\label{varvalp}
\end{table}

%%%%%%%%%%%%%%%%%%%%%%%%%%%%%%%%%%%%%%%%%%%%%%%%%%%%%%%%%%%%%%%%%%%%%%%%%%%%%%%%%%%%%%%
\begin{table}[h!]
	%\centering
		\begin{tabular}{|c|c|c||c|c|}
		\hline
			& \multicolumn{2}{|c||}{n=100} & \multicolumn{2}{c|}{n=1000}\\
			\hline
			&SRSWR&stratified&SRSWR&stratified\\
			\hline
		 $\|Var(\widehat{v_1})\|$&0.450&0.286&0.0396&0.0265\\
			\hline
			  $\|AV(\widehat{v_1})\|$&0.3997&0.287&0.0386&0.0267\\
			\hline
			  $\|\widehat{V}_p \left(\widehat{v}_1\right)\|$ & [0.335;0.491]&[0.252;0.354]&[0.0371;0.0410]&	[0.0256;0.0280]\\
			 \hline
		\end{tabular}
		
	\caption{Norm of the variance approximation of the first eigenvector estimator.} 
	\label{varvectp}
\end{table}	
%%%%%%%%%%%%%%%%%%%%%%%%%%%%%%%%%%%%%%%%%%%%%%%%%%%%%%%%%%%%%%%%%%%%%%%%%%%
\clearpage
\section*{Appendix : proofs}

%%%%%%%%%%%%%%%%%%%%%%%%%%%%%%%%%%%%%%%%%%
\begin{proof} of proposition \ref{distmesure}. \\ 
Let us introduce $\alpha_k = \frac{I_k}{\pi_k}-1,$ we have 
$$\begin{array}{lcl}
\displaystyle \frac{\widehat{N}- N}{N} & = & \displaystyle \frac{1}{N} \sum_{k \in U} \alpha_k .
\end{array}
$$
Noting that with assumptions (A2) and (A3),
$E(\alpha_k^2) = (1-\pi_k)/\pi_k < (1 - \pi_k)/\lambda,$ $| E(\alpha_k \alpha_\ell) | = | \Delta_{k \ell}/(\pi_k \pi_\ell) | \leq | \Delta_{k \ell} | / \lambda^2,$
and  taking now the expectation, according to the sampling distribution $p,$ we get 
\begin{eqnarray}
E_p  \left(\frac{\widehat{N}- N}{N}  \right)^2 & = & \frac{1}{N^2}  \sum_{k,\ell \in U} E_p (\alpha_\ell \alpha_k) \nonumber  \\
 & = & \frac{1}{N^2} \left(  \sum_{k \in U} \frac{1- \pi_k}{\pi_k} + \sum_{k \in U} \sum_{\ell \neq k} \frac{\Delta_{k \ell}}{\pi_k \pi_\ell} \right)  \nonumber\\
 & \leq & \frac{1}{N^2} \left(  \frac{N}{\lambda} + \frac{N(N-1)}{n} \frac{n \max |\Delta_{k \ell}|}{\lambda^2} \right)  \nonumber \\
 & = & O \left(\frac{1}{n} \right)
\label{borneNchap}
\end{eqnarray}
which is the first result. Looking now at the estimator of the mean function, we have
$$\begin{array}{lcl}
\displaystyle \widehat{\mu}- \mu & = & \displaystyle \frac{1}{N} \sum_{k \in U} \alpha_k Y_k + \left( \frac{1}{\widehat{N}} - \frac{1}{N} \right) \sum_{k \in s} \frac{1}{\pi_k} Y_k \\
 & = & \displaystyle \frac{1}{N} \sum_{k \in U} \alpha_k Y_k + \left( \frac{N-\widehat{N}}{N} \right) \widehat{\mu}
\end{array}
$$
By assumptions (A1)-(A3) it is clear that $\| \widehat{\mu} \| =O(1)$ and consequently $E_p \left\| \frac{N-\widehat{N}}{N} \widehat{\mu}\right\|^2 = O(n^{-1}).$
The first term of the right side of the inequality is dealt with as in (\ref{borneNchap}), noticing that $\|Y_k\| \leq C$ for all $k$:
\begin{eqnarray*}
E_p \left\| \frac{1}{N} \sum_{k \in U} \alpha_k Y_k \right\|^2  & = &   \frac{1}{N^2} \sum_{k, \ell \in U} E_p \left( \alpha_\ell \alpha_k \right) \ \langle Y_k, Y_\ell \rangle  \nonumber \\
& \leq & \frac{1}{N^2} \sum_{k,\ell \in U}  \left| E_p (\alpha_\ell \alpha_k) \right| \left\| Y_k \right\|   \left\|Y_\ell \right\|  \nonumber \\
 & = & O \left(\frac{1}{n} \right).
\end{eqnarray*}
To complete the proof, let us introduce the operator  $Z_k = Y_k \otimes Y_k$ and remark that 
\begin{eqnarray*}
\widehat{\Gamma} - \Gamma & = & \frac{1}{N} \sum_{k \in U}  \alpha_k Z_k  +  \left( \frac{1}{\widehat{N}} - \frac{1}{N} \right) \sum_{k \in s} \frac{1}{\pi_k} Z_k  + \mu \otimes \mu - \widehat{\mu} \otimes \widehat{\mu}
\end{eqnarray*}
By assumption (A1), we have that $\left| \langle Z_k, Z_\ell \rangle_{2} \right| \leq \left\| Y_k \right\|^2  \left\| Y_\ell  \right\|^2\leq C^4,$ for all $k$ and $\ell$ and we get with similar arguments as above that 
$$
E_p \left\| \frac{1}{N} \sum_{k \in U}  \alpha_k Z_k \right\|_2^2= O \left( \frac{1}{n} \right)
$$ and $E_p  \left\|  \left( \frac{1}{\widehat{N}} - \frac{1}{N} \right) \sum_{k \in s} \frac{1}{\pi_k} Z_k \right\|_2^2 = O(n^{-1}).$
Remarking now that 
$$\left\|  \mu \otimes \mu - \widehat{\mu} \otimes \widehat{\mu} \right\|_{2} \leq \left\| (\mu - \widehat{\mu})\otimes \mu  \right\|_2 +  \left\| \widehat{\mu}\otimes \left(\mu - \widehat{\mu} \right) \right\|_2
$$ 
the result is proved.

Consistency of the eigenelements is an immediate consequence of classical properties of the eigenelements of covariance operators. The  eigenvalues (see \textit{e.g.} Dauxois \textit{et al.}, 1982) satisfy
$| \widehat{\lambda}_j - \lambda_j | \leq \left\| \widehat{\Gamma} - \Gamma  \right\|_{2} \ .$
On the other hand, Lemma 4.3 by Bosq (2000) tells us that
$\left\| \widehat{v}_j - v_j \right\|  \leq C \delta_j \left\| \widehat{\Gamma} - \Gamma  \right\|_{2}$ 
where $\delta_1 = 2 \sqrt{2} (\lambda_1 - \lambda_2)^{-1}$  
and for $j \geq 2,$ 
\begin{eqnarray}
\delta_j &= & \displaystyle 2 \sqrt{2} \max \left[(\lambda_{j-1} - \lambda_j)^{-1},(\lambda_j - \lambda_{j+1})^{-1} \right].
\label{def:deltaj}
\end{eqnarray}
This concludes the proof.
\end{proof}

\begin{proof} of proposition \ref{IFacpf} : \\
Considering first the mean  curve $\mu,$ we get directly 
\begin{eqnarray*}
\mu(M + \epsilon \delta y) & = & \frac{1}{N+ \epsilon} \left( \sum_{\ell \in U} Y_\ell + \epsilon y \right) =  \mu + \frac{\epsilon}{N} \left( y - \mu \right) + o(\epsilon),
\label{def:FImu}
\end{eqnarray*}
so that
\begin{eqnarray*}
I \mu (M, Y_k)  & =  &\frac{1}{N}(Y_k - \mu).
\end{eqnarray*}
Let us first note that perturbation theory (Kato, 1966, Chatelin 1983) allows us to get the influence function of the eigenelements provided the influence function of the covariance operator is known. Indeed, let us consider the following expansion of $\Gamma$ according to some operator $\Gamma_1,$
\begin{eqnarray}
\Gamma(\epsilon) & = & \Gamma + \epsilon \Gamma_1 + o(\epsilon),
\label{devl:gamma}
\end{eqnarray}
we get from perturbation theory that the eigenvalues satisfy \begin{eqnarray}
\lambda_j(\epsilon) &=& \lambda_j + \epsilon \ \mbox{tr}\left(\Gamma_1 P_j\right)  + o(\epsilon),
\label{Kato:lambda}
\end{eqnarray}
where $P_j = v_j \otimes v_j$ is the projection onto the space spanned by $v_j$ and the trace of an operator $\Delta$ defined on $L^2[0,1]$ is defined by tr$(\Delta) = \sum_j  \langle \Delta e_j,e_j \rangle $ for any orthonormal  basis $e_j, j \geq 1$ of $L^2[0,1].$  There exists a similar result for the eigenfunctions which states, provided $\epsilon$ is small enough and for simplicity that the non null eigenvalues are distinct, that
\begin{eqnarray}
v_j(\epsilon) &=& v_j + \epsilon \left( S_j \Gamma_1(v_j) \right) + o(\epsilon),
\label{Kato:vecp}
\end{eqnarray}
where operator $S_j$ is defined on $L^2[0,1]$ as follows $$
S_j \ = \  \sum_{\ell \neq j}  \frac{ v_\ell \otimes v_\ell }{ \lambda_j - \lambda_\ell} \ .
$$
So going back to the notion of influence function, if we get an expression for  $\Gamma_1$ in our case, we will be able to derive the influence function for the eigenelements.
The influence function of $\Gamma$ can be computed directly using the definition, 
\begin{eqnarray}
\Gamma(\epsilon) & =  &\Gamma(M+ \epsilon \delta_y)  \nonumber \\
 & = & \frac{1}{N+\epsilon}  \left( \sum_{\ell \in U} \left( Y_\ell \otimes Y_\ell \right)  + \epsilon ( y  \otimes y)\right)  - \frac{1}{(N+\epsilon)^2} (N \mu + \epsilon y) \otimes (N \mu + \epsilon y) \nonumber \\
  & = & \Gamma + \frac{\epsilon}{N} \left( y \otimes y - \mu \otimes \mu - \Gamma \right) - \frac{\epsilon}{N} \left( \mu \otimes (y - \mu) + (y - \mu) \otimes \mu  \right) \nonumber  + o(\epsilon)\\
   & = & \Gamma + \frac{\epsilon}{N}  \left( (y - \mu) \otimes( y - \mu) - \Gamma \right) + o(\epsilon)
\label{def:FIGamma}
\end{eqnarray}
so that
\begin{eqnarray}
 I \Gamma (M,Y_k) & = & \frac{1}{N} \left( (Y_k-\mu)\otimes (Y_k-\mu) - \Gamma \right).\nonumber
\end{eqnarray}

The combination of (\ref{Kato:lambda}) and (\ref{def:FIGamma}) give us the influence function
of the $j$th eigenvalue 
\begin{eqnarray*}
I \lambda_j (M,Y_k) &= & \frac{1}{N} \left( \langle Y_k - \mu, v_j \rangle^2 - \lambda_j \right)
\end{eqnarray*}
as well as the influence function of the $j$th eigenfunction (since $\langle v_j , v_\ell \rangle = 0$ when $j \neq \ell$)
\begin{eqnarray*}
I v_j (M,Y_k) &= & \frac{1}{N} \left( \sum_{\ell \neq j} \frac{\langle Y_k - \mu, v_j \rangle \langle Y_k - \mu, v_\ell \rangle}{\lambda_j - \lambda_\ell} v_\ell\right).
\end{eqnarray*}
\end{proof}

%%%% Developpement a la von Mises
\begin{proof}  of proposition \ref{Varapprox}. \\
Let us begin with the mean function. The remainder term is defined as follows 
$$R_\mu = \widehat{\mu} - \mu - \int  I\mu (M, Y) d(\widehat{M}-M)$$ 
and 
$$ \begin{array}{lcl}
R_\mu& = & \displaystyle \widehat{\mu} - \mu -   \frac{1}{N} \sum_{k \in s} \frac{Y_k - \mu}{\pi_k} \\
 & = & \widehat{\mu} \left( 1 - \frac{ \widehat{N}}{N} \right) + \mu \left(  \frac{\widehat{N}}{N}  - 1 \right) \\
 & = & \left( \mu -  \widehat{\mu} \right) \left(  \frac{\widehat{N}}{N}  - 1 \right) \\
 &= & o_p(n^{-1/2}),
\end{array}
$$
since $\mu -  \widehat{\mu} = O_P(n^{-1/2})$ and $(\widehat{N} - N)/N = O_P(n^{-1/2})$ by proposition \ref{distmesure}.

For the covariance operator, we have
\begin{eqnarray}
R_\Gamma & = & \displaystyle \widehat{\Gamma} - \Gamma -  \frac{1}{N} \sum_{k \in s} \frac{1}{\pi_k} \left( (Y_k - \mu) \otimes (Y_k - \mu) - \Gamma \right)  \nonumber\\
 &= & \displaystyle \Gamma \left( \frac{\widehat{N}}{N}  - 1 \right) + \widehat{\Gamma} -  \frac{1}{N} \sum_{k \in s} \frac{1}{\pi_k} (Y_k - \mu) \otimes (Y_k - \mu) \nonumber \\
  & =  & \displaystyle \left( \Gamma - \widehat{\Gamma} \right) \left( \frac{\widehat{N}}{N}  - 1 \right) - \frac{\widehat{N}}{N} \left( (\mu - \widehat{\mu})  \otimes (\mu - \widehat{\mu}) \right)  \nonumber \\
 & = & o_p(n^{-1/2}), 
\label{ResteGamma}
\end{eqnarray}
noticing that 
$$
\frac{1}{N} \sum_{k \in s} \frac{ Y_k \otimes Y_k}{\pi_k} = \frac{\widehat{N}}{N} \left( \widehat{\Gamma} + \widehat{\mu} \otimes \widehat{\mu} \right).
$$
To study the remainder terms for the eigenelements, we need to go back to the perturbation theory and equations (\ref{devl:gamma}), (\ref{Kato:lambda})
and (\ref{Kato:vecp}). According to (\ref{ResteGamma}), with $\epsilon=n^{-1/2},$
we can write
\begin{eqnarray}
\Gamma_1  & = & \sqrt{n} \left( \frac{1}{N} \sum_{k \in s} \frac{1}{\pi_k} \left( (Y_k - \mu) \otimes (Y_k - \mu) - \Gamma \right) + R_\Gamma \right). 
\label{Gamma1:reste}
\end{eqnarray}
Introducing now (\ref{Gamma1:reste}) in equation (\ref{Kato:lambda}), we get noting that $\langle R_\Gamma  v_j, v_j \rangle = o_p(n^{-1/2}),$
\begin{eqnarray*}
\widehat{\lambda}_j - \lambda_j  & = & \frac{1}{N} \sum_{k \in s} \frac{1}{\pi_k} \left( \langle Y_k - \mu, v_j \rangle^2 -  \langle \Gamma v_j, v_j \rangle\right) + o_p(n^{-1/2})  \nonumber \\
 & = & \int  I\lambda_j (M, Y) d(\widehat{M}-M) + o_p(n^{-1/2})
\end{eqnarray*}
which proves that $R_{\lambda_j} = o_p(n^{-1/2}).$
Using now (\ref{Kato:vecp}) and since $S_j R_\Gamma  v_j = o_p(n^{-1/2})$, we can check with similar arguments that
\begin{eqnarray*}
\widehat{v}_j - v_j & = &  S_j   \left( \frac{1}{N} \sum_{k \in s} \frac{1}{\pi_k} \left( \langle Y_k - \mu, v_j \rangle (Y_k - \mu) - \lambda_j v_j \right) \right)  + o_p(n^{-1/2}) \nonumber \\
 & = &  \frac{1}{N} \sum_{k \in s}  \frac{1}{\pi_k} \sum_{\ell \neq j} \frac{ \langle Y_k - \mu, v_j \rangle \langle Y_k - \mu, v_\ell \rangle}{\lambda_j - \lambda_\ell}  v_\ell + o_p(n^{-1/2}) \nonumber \\
 &=& \int  Iv_j (M, Y) d(\widehat{M}-M) + o_p(n^{-1/2})
\end{eqnarray*}
and the proof is complete.
\end{proof}

\begin{proof} of proposition \ref{estimvarconv}. \\

We prove the result for functional linearized variables $u_k.$ For real valued linearized variables, for instance for an eigenvalue $\lambda_j$, the proof is similar  replacing the tensor product with  usual product and the norm $||\cdot||_{2}$ with the absolue value $|\cdot|.$
Let us denote by 
\begin{eqnarray*}
\widehat{AV}(T(\widehat M))=\sum_s\sum_s\frac{\Delta_{kl}}{\pi_{kl}}\frac{u_k}{\pi_k}\otimes\frac{u_l}{\pi_l}=\sum_U\sum_U\frac{\Delta_{kl}}{\pi_{kl}}\frac{u_k}{\pi_k}\otimes\frac{u_l}{\pi_l}I_kI_l
\end{eqnarray*}
and by
\begin{eqnarray*}
A &= &  \left\| AV(T(\widehat M))-\widehat{AV}(T(\widehat M))\right\|_{2} \quad \mbox{and} \quad B=  \left\| \widehat{AV}(T(\widehat M))-\widehat{V}_p(T(\widehat M))\right\|_{2}.
\end{eqnarray*}
It is clear that
\begin{eqnarray*}
\left\|AV(T(\widehat M))-\widehat{V}_p(T(\widehat M))\right\|_{2} & \leq & A+B.
\end{eqnarray*}
Let us consider
$$
E_p \left(A^2\right)= 
%E_p \left\langle \sum_{ k \in U}\sum_{ l \in U} \Delta_{kl} \left( 1 - \frac{I_k I_l}{\pi_{kl}} \right) \frac{u_k}{\pi_k}\otimes\frac{u_l}{\pi_l}, \sum_{ k \in U}\sum_{ l \in U} \Delta_{kl} \left( 1 - \frac{I_k I_l}{\pi_{kl}} \right) \frac{u_k}{\pi_k}\otimes\frac{u_l}{\pi_l} \right\rangle_2  
\sum_{ k,l \in U}\sum_{k',l' \in U}   \Delta_{kl}\Delta_{k'l'}E_p  \left( 1 - \frac{I_k I_l}{\pi_{kl}} \right) \left( 1 - \frac{I_{k'} I_{l'}}{\pi_{k'l'}} \right)   \left\langle  \frac{u_k}{\pi_k}\otimes\frac{u_l}{\pi_l}, \frac{u_{k'}}{\pi_{k'}}\otimes\frac{u_{l'}}{\pi_{l'}} \right\rangle_2 .
% & \leq &   E_p \left| \Delta_{kl} \left( 1 - \frac{I_k I_l}{\pi_{kl}} \right)  \right|^2 \left\|\frac{u_k}{\pi_k}\otimes\frac{u_l}{\pi_l} \right\|_2^2 + 2 \sum_{\ell \in U} \sum_{k,l , k \neq l} \\
 %& \leq & \sum_{ k \in U}\sum_{ l \in U} \Delta_{kl} \left( E_p \left( 1 - \frac{I_k I_l}{\pi_{kl}} \right)^2 \right)^{1/2}\left\|\frac{u_k}{\pi_k}\otimes\frac{u_l}{\pi_l} \right\|_2 \\
% & \leq & \left( \sum_{ k \in U}\sum_{ l \in U} \Delta_{kl}^2 E_p \left( 1 - \frac{I_k I_l}{\pi_{kl}} \right)^2 \right)^{1/2} \left(  \sum_{ k \in U}\sum_{l \in U} \left\|\frac{u_k}{\pi_k}\otimes\frac{u_l}{\pi_l} \right\|_2^2 \right)^{1/2} \\
%  & \leq & A_1 \ A_2.
$$
Using the fact that $\| u_k \otimes u_l \|_2 \leq \|u_k\|  \| u_l \|$ and since it is easy to check that $\|u_k\| <CN^{-1}$ where $C$ is a constant that does not depends on $k,$ we get, under assumptions  (A2), (A3) and (A4), with a similar decomposition as in Breidt and Opsomer (2000, proof of Th. 3) that $E_p \left(A^2\right) = o(n^{-2})$ and thus $E_p (A) = o(n^{-1}).$

Let us study now the second term $B$ and examine separately the case of the mean function and the eigenvalues and the case of the eigenfunctions which can not be dealt with the same way.
We can prove, under assumptions (A2) and (A3),  with similar manipulations as before that there exist some positive constant $C_2, C_3, C_4$ and $C_5$ such that
\begin{eqnarray*}
E_p \left(B \right) & =  & E_p \left\|  \sum_{k \in U}  \sum_{l \in U} \frac{\Delta_{kl}}{\pi_{kl}} I_kI_l \left( \frac{u_k}{\pi_k}\otimes\frac{u_l}{\pi_l} - \frac{\hat u_k}{\pi_k}\otimes\frac{\hat u_l}{\pi_l} \right) 
 \right\|_2 \\
 & \leq & \sum_{k \in U}  \sum_{l \in U} E_p \left| \frac{\Delta_{kl}}{\pi_{kl}} \right| I_k I_l \left\|  \frac{u_k}{\pi_k}\otimes\frac{u_l}{\pi_l} - \frac{\hat u_k}{\pi_k}\otimes\frac{\hat u_l}{\pi_l} \right\|_2 \\
& \leq & \sum_{k \in U}  \sum_{l \in U} \left( E_p \left(\frac{\Delta_{kl}}{\pi_{kl}}  I_k I_l  \right)^2 \right)^{1/2} \left( E_p \left\|  \frac{u_k}{\pi_k}\otimes\frac{u_l}{\pi_l} - \frac{\hat u_k}{\pi_k}\otimes\frac{\hat u_l}{\pi_l} \right\|_2^2 \right)^{1/2} \\ 
 & \leq & \frac{C_2}{N} \sum_{k \in U}  \sum_{l \neq k}  \left( E_p \left( \left\|  (u_k - \widehat{u}_k)\otimes u_l -\widehat u_k\otimes (\widehat u_l - u_l)\right\|_2^2 \right) \right)^{1/2} \\
  &  & + C_3 \sum_{k \in U}\left( E_p \left( \left\|  (u_k - \widehat{u}_k)\otimes u_k -\widehat u_k\otimes (\widehat u_k - u_k)\right\|_2^2 \right) \right)^{1/2} \\
 & \leq & \frac{C_4}{N} \sum_{k \in U}  \sum_{l \neq k }  \left(E_p  \left\|  u_k - \widehat{u}_k \right\|^2 \left\|u_l\right\|^2 + E_p \left\|\widehat u_l - u_l\right\|^2 \left\| \widehat u_k \right\|^2  \right)^{1/2} \\
 &  & + C_5 \sum_{k \in U}\left( E_p \left(  \left\|  u_k - \widehat{u}_k \right\|^2 (\left\|u_k\right\|^2 + \left\|\widehat{u}_k \right\|^2 \right) \right)^{1/2}
\end{eqnarray*}
For $k \neq l$ we have with assumption (A3) that $\Delta_{kl}^2 \leq CN^{-2}.$ Furthermore, since $n \leq \widehat{N} \leq n/\lambda,$
  the estimated linearized variables for the mean function satisfy $\left\| \widehat{u}_k \right\|^2 = O(n^{-2})$ uniformly in $k$ as well as for the   eigenvalues $(\widehat{u}_k )^2 = O(n^{-2}).$

For the mean function $\mu$ we have 
$$
u_k - \widehat{u}_k =  \frac{1}{N} (\widehat \mu - \mu) + \frac{1}{\widehat N} \frac{ \widehat N - N}{N} (Y_k - \widehat{\mu})
$$ 
and thus we easily get that  $E_p \|  u_k - \widehat{u}_k \|^2 = O(N^{-3})$ uniformly in $k.$
Considering the eigenvalues, we have $$
u_k - \widehat{u}_k = \frac{1}{N} (\langle Y_k - \mu,v_j \rangle^2 -  \langle Y_k - \hat \mu, \hat v_j \rangle^2 + \widehat{\lambda}_j - \lambda_j) - \frac{1}{\widehat N} \frac{ \widehat N - N}{N} (\langle Y_k - \hat \mu, \hat v_j \rangle^2 - \widehat{\lambda}_j).
$$ 
After some manipulations we also get that $E_p (u_k - \widehat{u}_k)^2 = O(N^{-3})$ uniformly in $k.$ Combining the previous results we get $E_p (B) = o(n^{-1})$ and the result is proved.

The technique is different for the eigenfunctions $\widehat{v}_1, \ldots, \widehat{v}_q$ because we cannot bound easily terms like $E_p (\widehat\lambda_j -\widehat\lambda_{j+1})^{-1}$ which appear in the estimators of the linearized variables. By the Cauchy Schwarz inequality we have  
\begin{eqnarray}
B & \leq & \left( \sum_{k \in U} \sum_{l \neq k } \left( \frac{ \Delta_{kl} I_k I_l}{\pi_{kl} \pi_k \pi_l} \right)^2 \right)^{1/2}
 \left( \sum_{k \in U} \sum_{l \neq k } \left\| u_k \otimes u_l -  \widehat{u}_k \otimes \widehat{u}_l \right\|^2_2 \right)^{1/2}   \label{cauchy:B1}
\\
  &  & + \left( \sum_{k \in U}  \left( \frac{ \Delta_{kl} I_k }{\pi_{kl} \pi_k^2 } \right)^2 \right)^{1/2}
 \left( \sum_{k \in U}  \left\| u_k \otimes u_k -  \widehat{u}_k \otimes \widehat{u}_k \right\|^2_2 \right)^{1/2}  .
 \label{cauchy:B2}
 \end{eqnarray} 
 By assumptions (A2) and (A3) we have, for $k \neq l,$
 $$
 E_p \left( \frac{ \Delta_{kl} I_k I_l}{\pi_{kl} \pi_k \pi_l} \right)^2 \ = \ \frac{\Delta_{kl}^2}{\pi_{kl} \pi_k^2 \pi_l^2} 
 \ \leq \  \frac{C_6}{n^2},
 $$
for some constant $C_6$ that does not depend on $k$ and $l.$
When $k=l,$ we have $ E_p \left( \frac{ \Delta_{kk} I_k}{ \pi_k^3} \right)^2 \leq C_7.$ Thus, by Markov inequality we have
\begin{eqnarray*}
\left( \sum_{k \in U} \sum_{l \neq k} \left( \frac{ \Delta_{kl} I_k I_l}{\pi_{kl} \pi_k \pi_l} \right)^2 \right)^{1/2} & = &  O_p(1),
\end{eqnarray*}
and 
\begin{eqnarray*}
\left( \sum_{k \in U}  \left( \frac{ \Delta_{kk} I_k }{\pi_{k}^3} \right)^2 \right)^{1/2} & = &  O_p(\sqrt{n}).
\end{eqnarray*}

Considering the terms containing linearized variables  in (\ref{cauchy:B1}) and (\ref{cauchy:B2}), we have the general inequality
\begin{eqnarray*}
\sum_{k \in U} \sum_{l \neq k} \left\| u_k \otimes u_l -  \widehat{u}_k \otimes \widehat{u}_l \right\|^2_2  & \leq & 2 \sum_{k \in U}  \sum_{l \neq k}  \  \left\|  u_k - \widehat{u}_k \right\|^2 \left\|u_l\right\|^2 + \left\|\widehat u_l - u_l\right\|^2 \left\| \widehat u_k \right\|^2 .
\end{eqnarray*}
Let us make now  the following decomposition
\begin{eqnarray}
\left\| u_k - \widehat{u}_k \right\| & \leq & \left\| N u_k\right\| \left( \frac{\widehat N - N}{N \widehat N} \right) + \frac{1}{\widehat N} \left\| \widehat N \widehat u_k - N u_k \right\| 
\label{decukvk}
\end{eqnarray}
with
$$
N u_k - \widehat N \widehat{u}_k = \langle Y_k - \mu, v_j \rangle \sum_{\ell \neq j} \frac{\langle Y_k - \mu, v_\ell \rangle}{\lambda_j - \lambda_\ell} v_\ell   - \langle Y_k - \hat \mu, \hat v_j \rangle \sum_{\ell \neq j}  \frac{ \langle Y_k - \hat \mu, \hat v_\ell \rangle}{\hat \lambda_j - \hat \lambda_\ell} \hat v_\ell .
$$
It is clear that $\left\| N u_k\right\|=O(1)$ uniformly in $k$ and $\left( \frac{\widehat N - N}{N \widehat N} \right)=O_p(n^{-3/2}).$
We have for the second right hand term of  inequality (\ref{decukvk}),
\begin{eqnarray}
 \left\| N u_k - \widehat N \hat u_k \right\| &\leq & \left|Ê\langle Y_k - \mu, v_j \rangle - \langle Y_k - \hat \mu, \hat v_j \rangle\right| \left\| \sum_{\ell \neq j} \frac{\langle Y_k - \mu, v_\ell \rangle}{\lambda_j - \lambda_\ell} v_\ell \right\| \nonumber\\
 &  + &  \left| \langle Y_k - \hat \mu, \hat v_j \rangle \right| \left\| \sum_{\ell \neq j} \frac{\langle Y_k - \mu, v_\ell \rangle}{\lambda_j - \lambda_\ell} v_\ell - \sum_{\ell \neq j} \frac{ \langle Y_k - \hat \mu, \hat v_\ell \rangle}{\hat \lambda_j - \hat \lambda_\ell} \hat v_\ell \right\| . \label{Nuk2}
\end{eqnarray}
It is clear that the first term at the right hand side of previous inequality satisfies, uniformly in $k,$ 
$$
E_p \left( \left| \langle Y_k - \mu, v_j \rangle - \langle Y_k - \hat \mu, \hat v_j \rangle \right| \ \left\| \sum_{\ell \neq j} \frac{\langle Y_k - \mu, v_\ell \rangle}{\lambda_j - \lambda_\ell} v_\ell \right\| \right)^2 = O \left( \frac{1}{n} \right).
$$
Let us introduce the random variable 
$$T = \min(\lambda_j - \lambda_{j+1}, \lambda_{j-1} - \lambda_j) \min(\hat \lambda_j - \hat \lambda_{j+1}, \hat \lambda_{j-1} - \hat \lambda_j ),$$
the eigenvalues being distinct, we have  with Proposition 1 that
$
\frac{1}{T} = O_p(1).
$
 As far as the second term in (\ref{Nuk2}) is concerned we can write
\begin{eqnarray}
& & \left\| \sum_{\ell \neq j} \frac{(\hat \lambda_j - \hat \lambda_\ell)\langle Y_k - \mu, v_\ell \rangle v_\ell - (\lambda_j - \lambda_\ell)  \langle Y_k - \hat \mu, \hat v_\ell \rangle  \hat v_\ell}{(\lambda_j - \lambda_\ell)(\hat \lambda_j - \hat \lambda_\ell)} \right\|^2 \nonumber\\
& \leq & \frac{4 \left| \lambda_j  - \hat \lambda_j \right|^2 }{T^2}    \sum_{\ell \neq j}  \langle Y_k - \mu, v_\ell \rangle^2 + \frac{4}{T^2}  \sum_{\ell \neq j} (\lambda_\ell  - \hat \lambda_\ell)^2   \langle Y_k - \mu, v_\ell \rangle^2 \nonumber \\
 &  & + 4 \lambda_j^2   \left\|  \sum_{\ell \neq j} \frac{ \langle Y_k - \mu, v_\ell \rangle v_\ell - \langle Y_k - \widehat \mu, \widehat v_\ell \rangle \widehat v_\ell }{(\lambda_j - \lambda_\ell)(\hat \lambda_j - \hat \lambda_\ell)} \right\|^2 .
\label{Nuk3}
\end{eqnarray}
We have seen  that $\sup_\ell | \lambda_\ell - \widehat \lambda_\ell |^2 \leq \left\| \Gamma - \widehat \Gamma \right\|^2$ and thus the first two terms in (\ref{Nuk3}) are $O_p(n^{-1}).$ 
The assumption that $\Gamma$ is a finite rank operator is needed to deal with the last term of (\ref{Nuk3}). 
Using the fact that $\| v_\ell - \widehat v_\ell \| \leq C \delta_j \left\| \Gamma - \widehat \Gamma \right\|$ where $\delta_j$ is defined in (\ref{def:deltaj}), we also get that this last term is also $O_p(n^{-1}).$
Combining all these results we finally get that, uniformly in $k$
$$
\left\| u_k - \widehat u_k \right\| = O_p(n^{-3/2}).
$$
It can be checked easily, under the finite rank assumption of $\Gamma$ that, uniformly in $k,$  $\| \widehat u_k \| =O_p(n^{-1})$ and $\| u_k \| =O(n^{-1}).$
Going back now to (\ref{cauchy:B1}) and (\ref{cauchy:B2}) we get that $B = O_p(1)O_p(n^{-3/2}) + O_p(n^{1/2})O_p( n^{-2}) = o_p(n^{-1}).$
This concludes the proof.
\end{proof}

\noindent \textbf{Acknowledgments.} We would like to thank Andr\'e Mas for helpful comments as well as the two referees for their constructive remarks that helped us to improve the manuscript.

\section*{References}
\begin{description}

 \item Benko, M., H\"ardle, W. and Kneip, A. (2009). Common functional principal components. {\it Annals of Statistics}, {\bf 37}, 1-34.

\item Berger, Y.G, Skinner, C.J (2005).
\newblock A jacknife variance estimator for unequal probability sampling.
\newblock {\em J. R. Statist. Soc B}, {\bf 67}, 79-89.

\item Besse, P.C and  Ramsay, J.O. (1986).
\newblock Principal component analysis of sampled curves.
\newblock {\em Psychometrika}, {\bf 51}, 285-311.

\item  Besse, P.C., Cardot, H. and Stephenson, D.B. (2000). 
\newblock Autoregressive Forecasting of Some Functional Climatic Variations. 
\newblock {\em Scand. J. Statist.}, \textbf{27}, 673-687.

\item Bosq, D. (2000). 
\newblock {\em Linear Processes in Function Spaces.} 
\newblock Lecture Notes in Statistics,  149, Springer. 

\item Breidt, F.J. and Opsomer, J.D. (2000). Local Polynomial Survey Regression Estimators in Survey Sampling. \textit{The Annals of Statistics}, {\bf 4}, 1026-1053.

\item Campbell, C. (1980). A Different View of Finite Population Estimation.
\textit{Proceeding of the Section on Survey Research Methods}, American Statistical Association. 319-324.

\item Cardot, H. (2000).
\newblock Nonparametric estimation of the smoothed principal components
  analysis of sampled noisy functions.
\newblock {\em J. Nonparametr. Stat.}, {\bf 12}, 503-538.

\item Cardot, H., Faivre, R. and  Goulard, M. (2003).
\newblock  Functional approaches for predicting land use with the temporal
     evolution of coarse resolution remote sensing data. 
\newblock {\em J. of Applied Statistics}, {\bf 30}, 1185-1199.

\item Cardot, H., Josserand, E. (2009). Sondages stratifi\'es pour donn\'ees fonctionnelles : allocation optimale et bandes de confiance asymptotiques. {\em Preprint.}

\item Castro, P., Lawton, W. and Sylvestre, E. (1986).
Principal Modes of Variation for Processes with Continuous Sample Curves.
{\em Technometrics}, {\bf 28}, 329-337.

\item Chatelin, F. (1983). \textit{Spectral approximation of linear operators.} Academic Press, New York

\item Chen, J., Rao, J.N.K. (2007). Asymptotic Normality Under Two-phase Sampling Designs. {\it Statistica Sinica}, {\bf 17}, 1047-1064.

\item Chiky, R., H\'ebrail, G. (2009). Spatio-temporal sampling of distributed data streams. {\em J. of Computing Science and Engineering}, to appear.

\item Chiou, J-M.,  M\"uller, H.G. and Wang, J.L. (2003). 
\newblock Functional  quasi-likelihood regression models with smooth random effects.
\newblock  {\em J. Roy.  Statist. Soc. Ser. B}, {\bf 65}, 405-423.

\item Croux, C., Ruiz-Gazen, A. (2005). High breakdown estimators for principal components : the projection-pursuit approach revisited. \textit{J. Multivariate Analysis}, \textbf{95}, 206-226.

\item Cuevas, A., Febrero, M. and  Fraiman, R. (2002). Linear functional regression: The case of fixed design and functional response. \textit{Canadian Journal of Statistics}, {\bf 30}, 285-300.

 \item Dauxois, J., Pousse, A., and Romain, Y. (1982). 
\newblock Asymptotic theory for the principal component analysis of a random 
  vector function: some applications to statistical inference. 
\newblock {\em J. Multivariate Anal.}, {\bf 12}, 136-154. 

\item Davison, A.C. and Hinkley, D.V. (1997).
\textit{Bootstrap Methods and Their Application.} Cambridge: Cambridge University Press.

\item Dessertaine, A. (2006).
\newblock Sondage et s\'eries temporelles: une application pour la pr\'evision de la consommation \'electrique.
\textit{38\`emes Journ\'ees de Statistique}, Clamart, Juin 2006.

\item Deville, J.C. (1974).
\newblock M\'ethodes statistiques et num\'eriques de l'analyse harmonique.
\newblock {\em Ann. Insee}, {\bf 15}, 3-104.

 \item Deville, J.C. (1999).
\newblock Variance estimation for complex statistics and estimators: linearization and residual techniques.
\newblock {\em Survey Methodology}, {\bf 25}, 193-203.

\item Ferraty, F. and Vieu, P. (2006). \textit{Nonparametric Functional Data Analysis, Theory and Applications.} Springer Series in Statistics, Springer, New-York.

\item Hampel, F. R. (1974). The influence curve and its role in robust statistics. \textit{J. Am. Statist. Ass.}, \textbf{69}, 383-393. 

\item Hastie, T. and Mallows, C. (1993). A discussion of ``A statistical view of some chemometrics regression tools'' by I.E. Frank and J.H. Friedman. \textit{Technometrics}, \textbf{35}, 140-143.

\item Isaki, C.T. and Fuller, W.A. (1982). Survey design under the regression superpopulation model. {\em J. Am. Statist. Ass.} \textbf{77}, 89-96.

\item James, G., Hastie, T., and Sugar, C. (2000).
\newblock Principal Component Models for Sparse Functional Data.
 \newblock \textit{Biometrika},   {\bf 87} , 587-602.

\item Kato, T. (1966). \textit{Perturbation theory for linear operators.} Springer Verlag, Berlin.

\item Kirkpatrick, M. and  Heckman, N. (1989). A quantitative genetic model for growth, shape, reaction norms and other infinite dimensional characters. 
\newblock \textit{J. Math. Biol.}, \textbf{27}, 429-450

\item Kneip, A. and Utikal, K.J. (2001). Inference for Density Families Using Functional Principal Component Analysis. \textit{J. Am. Statist. Ass.}, {\bf 96}, 519-542.

\item Mises, R., v (1947). On the asymptotic distribution of differentiable statistical functions. \textit{Ann. Math. Statist.}, 18, 309-348.

\item M\"uller, H.G. and Yao, F. (2008). 
\newblock Functional  additive models. {\em Preprint.}
 
 \item Ramsay, J. O. and Silverman, B.W. (2002). \emph{Applied Functional Data Analysis: Methods and Case Studies}. Springer-Verlag.

\item Ramsay, J. O. and Silverman, B.W. (2005). \emph{Functional Data Analysis}. Springer-Verlag, second edition.

\item Rice, J. (2004). Functional and Longitudinal Data Analysis. {\em Statistica Sinica}, {\bf 14}, 613-629.

\item Robinson, P.M. and S\"{a}rndal, C.E. (1983). Asymptotic properties of the generalized regression estimator in probability sampling. \textit{Sankhya : The Indian Journal of Statistics}, {\bf 45}, 240-248.

\item S\"{a}rndal, C.E.,  Swensson, B.  and J. Wretman, J. (1992). {\em Model Assisted Survey Sampling}. Springer-Verlag.

\item Serfling, R. (1980). \textit{Approximation Theorems of Mathematical Statistics}, John Wiley and Sons.

 \item Skinner, C.J, Holmes, D.J, Smith, T.M.F  (1986).
\newblock The Effect of Sample Design on Principal Components Analysis.
\newblock {\em J. Am. Statist. Ass.} {\bf 81}, 789-798. 

\end{description}

\end{document}